\documentstyle[12pt]{article}
\begin{document}
\newtheorem{theorem}      {Th\'eor\`eme}[section]
\newtheorem{theorem*}     {theorem}
\newtheorem{proposition}  [theorem]{Proposition}
\newtheorem{definition}   [theorem]{Definition}
\newtheorem{e-lemme}        [theorem]{Lemma}
\newtheorem{cor}   [theorem]{Corollaire}
\newtheorem{resultat}     [theorem]{R\'esultat}
\newtheorem{eexercice}    [theorem]{Exercice}
\newtheorem{rrem}    [theorem]{Remarque}
\newtheorem{pprobleme}    [theorem]{Probl\`eme}
\newtheorem{eexemple}     [theorem]{Exemple}
\newcommand{\preuve}      {\paragraph{Preuve}}
\newenvironment{probleme} {\begin{pprobleme}\rm}{\end{pprobleme}}
\newenvironment{remarque} {\begin{rremarque}\rm}{\end{rremarque}}
\newenvironment{exercice} {\begin{eexercice}\rm}{\end{eexercice}}
\newenvironment{exemple}  {\begin{eexemple}\rm}{\end{eexemple}}
%
% english
%
\newtheorem{e-theo}      [theorem]{Theorem}
\newtheorem{theo*}     [theorem]{Theorem}
\newtheorem{e-pro}  [theorem]{Proposition}
\newtheorem{e-def}   [theorem]{Definition}
\newtheorem{e-lem}        [theorem]{Lemma}
\newtheorem{e-cor}   [theorem]{Corollary}
\newtheorem{e-resultat}     [theorem]{Result}
\newtheorem{ex}    [theorem]{Exercise}
\newtheorem{e-rem}    [theorem]{Remark}
\newtheorem{prob}    [theorem]{Problem}
\newtheorem{example}     [theorem]{Example}
\newcommand{\proof}         {\paragraph{Proof~: }}
\newcommand{\hint}          {\paragraph{Hint}}
\newcommand{\heuristicproof}{\paragraph{heuristic proof}}
\newenvironment{e-probleme} {\begin{e-pprobleme}\rm}{\end{e-pprobleme}}
\newenvironment{e-remarque} {\begin{e-rremarque}\rm}{\end{e-rremarque}}
\newenvironment{e-exercice} {\begin{e-eexercice}\rm}{\end{e-eexercice}}
\newenvironment{e-exemple}  {\begin{e-eexemple}\rm}{\end{e-eexemple}}
\newcommand{\1}        {{\bf 1}}
\newcommand{\pp}       {{{\rm I\!\!\! P}}}
\newcommand{\qq}       {{{\rm I\!\!\! Q}}}
\newcommand{\B}        {{{\rm I\! B}}}
\newcommand{\cc}       {{{\rm I\!\!\! C}}}
\newcommand{\N}        {{{\rm I\! N}}}
\newcommand{\R}        {{{\rm I\! R}}}
\newcommand{\D}        {{{\rm I\! D}}}
\newcommand{\Z}       {{{\rm Z\!\! Z}}}
\newcommand{\C}        {{\bf C}}        % ensemble des nombres complexes
\newcommand{\rank}{\hbox{rank}}
\newcommand{\CC}{{\cal C}}
\def\Re {{\rm Re\,}}
\def\Im {{ \rm Im\,}}
\def\st {{ \rm st}}
\def\ind {{ \rm ind}}
\def\const {{ \rm const}}
\def\A{{\cal A}}
\def\bar{\overline}

%
% utilitaires
%
\newcommand{\dontforget}[1]
{{\mbox{}\\\noindent\rule{1cm}{2mm}\hfill don't forget : #1
\hfill\rule{1cm}{2mm}}\typeout{---------- don't forget : #1 ------------}}
\newcommand{\note}[2]
{ \noindent{\sf #1 \hfill \today}

\noindent\mbox{}\hrulefill\mbox{}
\begin{quote}\begin{quote}\sf #2\end{quote}\end{quote}
\noindent\mbox{}\hrulefill\mbox{}
\vspace{1cm}
}
\title{ Boundary value problems and  $J$-complex curves }
\author{ Alexandre Sukhov{*} and Alexander Tumanov{**}}
\date{}
\maketitle

{\small

* Universit\'e des Sciences et Technologies de Lille, Laboratoire
Paul Painlev\'e,
U.F.R. de
Math\'e-matique, 59655 Villeneuve d'Ascq, Cedex, France,
 sukhov@math.univ-lille1.fr

** University of Illinois, Department of Mathematics
1409 West Green Street, Urbana, IL 61801, USA, tumanov@math.uiuc.edu
}
\bigskip

Abstract. We give a solution to the problem of filling
by a Levi-flat hypersurface for a class of totally real tori
in $\cc^2$ equipped with a certain almost complex structure.

MSC: 32H02, 53C15.

Key words: almost complex structure, Levi-flat hypersurface,
$J$-complex disc.
\bigskip

\section{Introduction}

The problem of gluing complex discs to real submanifolds
in complex space has been a subject of extensive research.
In this paper we consider the problem for a class of totally
real tori in $\cc^2$ equipped with a certain almost complex
structure. We restrict to structures $J$ with the following
characteristic property: the lines parallel to one
coordinate axis are $J$-complex hypersurfaces.
In complex dimension 2 every almost complex structure
$J$ locally has this property; we impose this condition globally.
We prove (Theorem \ref{F-Sh}) that certain real tori
can be filled by Levi-flat hypersurfaces.
The result can be viewed as a solution to
a Riemann-Hilbert problem for a quasi-linear elliptic equation
in the plane with non-linear boundary conditions.
For the usual complex structure in $\cc^2$, our result
was obtained earlier by Forstneri\v c \cite{Fo} and
Schnirelman \cite{Shni}.

Let $(M,J)$ be an almost complex manifold.
We denote by $\D$ the unit disc in $\cc$.
We denote by $J_\st$ the standard complex structure of $\cc^n$;
the value of $n$ will be clear from context.
A continuous map $f:\bar\D \to M$ differentiable in $\D$
is called a {\it $J$-complex} (or $J$-holomorphic) disc
if it satisfies the equation
$df \circ J_\st = J \circ df$ in $\D$.
We often identify a $J$-complex disc $f$ with the image
$f(\bar\D)$ and call it just a disc. By the boundary of
the disc $f$ we mean the restriction $f|_{b\D}$ which we also
identify with its image.

In local coordinates $Z =(z,w)\in\cc^2$, an almost complex
structure $J$ can be defined by a complex matrix $A(Z)$
so that a map $Z:\D \to U$ is $J$-complex if and only if it
satisfies the equation
\begin{eqnarray}
\label{CR}
Z_{\overline\zeta} - A(Z)\overline{Z_\zeta}=0.
\end{eqnarray}
The matrix $A(Z)$ is defined by
\begin{eqnarray}
\label{matrixA} A(Z)v = (J_\st + J(Z))^{-1}(J_\st -
J(Z))\overline v.
\end{eqnarray}
Indeed, it is easy to see that the right-hand side
of (\ref{matrixA}) is $\cc$-linear in $v\in\cc^2$,
hence it defines a unique matrix $A(Z)$
(see \cite{SuTu1} for details).
We call the matrix $A$ {\it the complex matrix} of $J$.
The ellipticity of (\ref{CR}) is equivalent to
$\det(I-A\bar A)\ne0$.
In a fixed coordinate chart, the  correspondence
between almost complex structures $J$ with $\det(J_\st+J)\ne0$
and complex matrices with $\det(I-A\bar A)\ne0$ is one-to-one
\cite{SuTu1}.

Let $A$ be a lower triangular complex matrix function
in a domain in $\cc^2$. That is,
\begin{eqnarray}
\label{A-abc}
A= \left(
\begin{array}{cll}
a & & 0\\
b & & c
\end{array}
\right)
\end{eqnarray}
The matrix $A$ is the matrix of an almost complex
structure if and only if $|a|\ne1$ and $|c|\ne1$.
An almost complex structure $J$ with matrix $A$ of
the form (\ref{A-abc}) has the following characteristic
property: the lines $z = \const$ are $J$-complex curves,
that is, their tangent spaces are $J$-invariant.

Our main result is the following
\begin{e-theo}
\label{F-Sh}
Let $J$ be a smooth almost complex structure
on $\bar\D\times\cc$ with matrix (\ref{A-abc}), in which
$|a|<a_0$, $|c|<a_0$, and $0<a_0<1$ is a constant.
Suppose the map $z\mapsto(z,0)$ is $J$-complex, that is,
$a(z,0)=b(z,0)=0$.
Suppose the first derivatives of $a$, $b$, and $c$
with respect to
$w$ and $\bar w$ are uniformly bounded.
Let $\gamma_z\subset\cc$ be a smooth simple closed curve
depending smoothly on the parameter $z\in b\D$.
Suppose that for every $z\in b\D$, the bounded
component of $\cc\setminus\gamma_z$ contains 0.
Introduce the torus $\Lambda=\bigcup_{z\in b\D}\{z\}\times\gamma_z$.
Then the following hold.
\begin{itemize}
\item[(i)]
For every point $p\in\Lambda$ there exists a  $J$-complex disc
$f:\bar\D\to\bar\D\times\cc$
such that $f(1)=p$ and whose image $f(\bar\D)$
is a graph of a non-vanishing function
$z\mapsto w(z)$ of class $C^{\infty}(\bar\D)$,
that is, $f(\bar\D)=\{(z,w): w=w(z), z\in\bar\D\}$.
In particular, $f$ is an embedding, $f(b\D)\subset \Lambda$,
and $f(\bar\D)$ does not meet $\bar\D\times\{0\}$.
The disc $f$ is unique up to parametrization
in any Sobolev class $W^{1,p}(\D)$, $p>2$.

\item[(ii)]
The discs in (i) form a $C^\infty$-smooth one-parameter family;
they depend continuously on the matrix $A$ and
the curves $\gamma_z$.
They are disjoint and fill a smooth Levi-flat
hypersurface $\Gamma \subset \bar\D\times\cc$
with boundary $\Lambda$.

\item[(iii)]
There is a constant $C>0$ that depends only on $A$
such that if $\gamma_z\subset r\bar\D$ for some $r>0$ and
all $z\in b\D$, then $\Gamma\subset\bar\D\times R\bar\D$, $R=Cr$.
\end{itemize}
\end{e-theo}

The hypothesis that the matrix $A$ has a special form
(\ref{A-abc})
is natural because it is tied to the special form
of the domain $\bar\D\times\cc$. In particular, it guarantees
that the torus $\Lambda$ is totally real. For a general
matrix $A$ the conclusion of the theorem fails.

As we mention above, in the integrable case ($A=0$),
this result was obtained earlier by Forstneri\v c \cite{Fo}
and Schnirelman \cite{Shni}.
Already in this special case, Theorem \ref{F-Sh} has important
connections. Berndtsson and Ransford \cite{BeRa} used a suitable
version in their proof of the corona theorem.
Slodkowski \cite{Sl} used similar ideas in his $\lambda$-lemma,
which has important applications in complex dynamics and
quasi-conformal maps. Recently Duval and Gayet \cite{DuGa}
obtained a new result on gluing holomorphic discs and annuli
to certain totally real tori in $\cc^2$. In the integrable case,
Theorem \ref{F-Sh} can be used to study envelopes of holomorphy
following Bedford and Klingenberg \cite{BeKl} and
Forstneri\v c \cite{Fo}.
Another direction concerns applications in symplecic and
contact geometry in the spirit of works of Gromov \cite{Gr},
Eliashberg \cite{El}, and others. Here Theorem \ref{F-Sh}
may be used in full generality, that is, for non-integrable
almost complex structures.
In a forthcoming paper we will use Theorem \ref{F-Sh}
for describing deformations of $J$-complex discs.
Finally, a situation related to
that in Theorem \ref{F-Sh} often arises in the theory
of holomorphic foliations, see for instance \cite{Br}.
We hope that our result will find further applications in
these directions and will be a useful tool in their development.

\section{Elliptic  estimates and the maximum principle}

We first consider a non-homogeneous Beltrami equation
\begin{eqnarray}
\label{Beltrami-v}
v_{\bar z}=q v_z + Q.
\end{eqnarray}
The following result must be well known, however
we could not find a precise reference.
For completeness we include a proof.
\begin{e-pro}
\label{Beltrami-Re0}
\begin{itemize}
\item[(i)]
Let $q$ and $Q$ be bounded functions in $\D$,
$|q|\le q_0<1$, $|Q|\le Q_0$, here
$q_0$ and $Q_0$ are constants. There exists $p > 2$ and a unique solution
$v$ of (\ref{Beltrami-v}) in the Sobolev class $W^{1,p}(\D)$
 with boundary conditions
$\Re v|_{b\D}=0$, $v(1)=0$.
Furthermore  $||v||_{C^\alpha(\bar\D)}\le C$ with
$\alpha = \frac{p-2}{p}$. Here $C>0$ and $p$, hence $\alpha$,
depend on $q_0$ and $Q_0$ only.

\item[(ii)] If in addition
$||q||_{C^{k,\beta}(\bar\D)}
+||Q||_{C^{k,\beta}(\bar\D)}\le Q_0$
for some $0<\beta<1$ and $k\ge 0$, then
$||v||_{C^{k+1,\beta}(\bar\D)}\le C$; here
$C>0$ depends on $\beta$, $k$, $q_0$ and
$Q_0$ only.
\end{itemize}
\end{e-pro}

The statement (ii) is not needed in this paper. We include it
for completeness and future references.

\proof
(i)
We use a classical method for solving the Beltrami
equation based on the Cauchy-Green operator $T$
and its modification
(\cite{Ve}, Theorem 3.29;
see also \cite{CoSuTu}, equation (5)):
\begin{eqnarray*}
Tu(z)=\frac{1}{2\pi i}\int_{\D}\frac{u(\zeta)
d\zeta\wedge d\bar\zeta}{\zeta-z}, \qquad
T_1u(z)=Tu(z)-\bar{Tu(z^*)}-2i\Im Tu(1),
\end{eqnarray*}
here $z^*:=1/\bar z$.
Then $v=T_1u$ solves the boundary value problem:
$v_{\bar z}=u$, $\Re v|_{b\D}=0$, $v(1)=0$.
We also use $S_1u(z)=\partial_z T_1u(z)$.
The operator $S_1$ is an isometry of $L^2(\D)$.

Put $u:=v_{\bar z}$. Then we have
$v=T_1 u$ and $v_z=S_1 u$. Then (\ref{Beltrami-v})
turns into
\begin{eqnarray}
\label{Beltrami-int-eq}
u=qS_1u+Q.
\end{eqnarray}
Since $|q|\le q_0<1$ and $||S_1||_{L^2}=1$,
then for $p>2$ close to 2, we have $||qS_1||_{L^p}<1$.
Hence (\ref{Beltrami-int-eq}) has a unique solution
$u\in L^p(\D)$. Then $v=T_1 u\in W^{1,p}(\D)$. 
Then $v\in C^\alpha$,
$\alpha=\frac{p-2}{p}$. The solution is unique
in $W^{1,p}$ and its $C^\alpha$ norm is estimated
in terms of $q_0$ and $Q_0$ only.

(ii) Let $\xi:\bar\D\to\bar\D$ be a Beltrami
homeomorphism satisfying the equation
$$
\xi_{\bar z}=q\xi_z, \qquad
\xi(1)=1.
$$
Since $q\in C^{k,\beta}$, then $\xi\in C^{k+1,\beta}$
is a diffeomorphism (see \cite{Ve}, Theorem 2.10).
By changing the independent variable to $\xi$,
the equation (\ref{Beltrami-v}) turns into
$$
v_{\bar\xi} \,\bar\xi_{\bar z}(1-|q|^2)=Q.
$$
Since $\xi\in C^{k+1,\beta}(\D)$, then the change of variable
$z\to\xi$ preserves the class $C^{k,\beta}(\D)$.
Hence
$$
v=T_1\frac{Q}{\bar\xi_{\bar z}(1-|q|^2)},
$$
here $T_1$ is applied with respect to the
variable $\xi$.
It shows that $v\in C^{k+1,\beta}(\D)$.
It is clear from the construction that the norm
of $v$ is estimated in terms of
$\beta$, $q_0$, and $Q_0$ only.
The proof is complete.
Q.E.D.
\medskip

We consider some properties of almost complex structures
with matrices of the form (\ref{A-abc}).
We begin with coordinate changes that preserve
this form.

\begin{e-pro}
\label{Change-abc}
Let $A$ be a complex matrix of a smooth almost complex
structure $J$ in a domain in $\cc^2$ with coordinates $Z=(z,w)$.
Suppose $A$ has the form (\ref{A-abc})
with $|a|<1$ and $|c|<1$.
Let $Z'=(z',w')$ be a smooth orientation preserving
coordinate change of the form $z'=z$, $w'=w'(z,w)$.
Then the matrix $A'$ of $J$ relative
to the new coordinates also has the form (\ref{A-abc}).
Furthermore, if
$w(z,0)=0$,
$a(z,0)=0$, and
$b(z,0)=0$,
then for the corresponding entries of $A'$ we have
$a'(z',0)=0$ and
$b'(z',0)=0$.
(Finally, if in addition
$w'_{\bar w}(z,0)=0$ and
$c(z,0)=0$,
then
$c'(z',0)=0$.)
\end{e-pro}
\proof
We represent the coordinate change in an equivalent form
$z=z'$, $w=w(z',w')$.
According to \cite{SuTu1}, the matrix $A$ changes to
\begin{eqnarray*}
A' = (Z_{Z'} - A \bar Z_{Z'})^{-1}
(A\bar Z _{\bar Z'} - Z_{\bar Z'}).
\end{eqnarray*}
Then the matrix $A'$ has the form (\ref{A-abc})
with entries
\begin{eqnarray*}
&& a'=a,\\
&& b'=(w_{w'}-c\bar w_{w'})^{-1}(b-a(w_{z'}-c\bar w_{z'})
+c\bar w_{\bar z'}-w_{\bar z'}),\\
&& c'=(w_{w'}-c\bar w_{w'})^{-1}(c\bar w_{\bar w'}-w_{\bar w'}).
\end{eqnarray*}
Note that $w_{w'}-c\bar w_{w'}\ne 0$ because the change
or variables is orientation preserving and $|c|<1$.
The conclusions now follow from the above
expressions of $a'$, $b'$, and $c'$. Q.E.D.
\medskip

For almost complex structures with matrix (\ref{A-abc}),
one can reduce the Cauchy-Riemann system to a single equation.

\begin{e-pro}
\label{CR1prop}
Let $\zeta\to(z(\zeta),w(\zeta))$
be a $J$-complex curve for an almost complex
structure with matrix (\ref{A-abc}),
where $|a|<1$ and $|c|<1$.
Suppose $z_\zeta \ne 0$.
Then the map locally can be represented by a
graph of a function $z\mapsto w(z)$ satisfying
the equation
\begin{eqnarray}
\label{CR1}
w_{\bar z}=a_1 w_z + c_1 \bar w_{\bar z} + b_1,
\end{eqnarray}
whose coefficients are determined by $A$.
The equation  is elliptic, in particular $|a_1|+|c_1|<1$.
Moreover, if $a(z,0)=0$ (resp. $b(z,0)=0$, $c(z,0)=0$),
then $a_1(z,0)=0$ (resp. $b_1(z,0)=0$, $c_1(z,0)=0$).
Conversely, if $z\mapsto w(z)$ satisfies (\ref{CR1}),
then its graph locally can be represented as
a parametrized $J$-complex curve.
Finally, the correspondence between the triples
$(a,b,c)$ with $|a|<1$ and $|c|<1$ and the triples
$(a_1,b_1,c_1)$ with $|a_1|+|c_1|<1$ is one-to-one.
\end{e-pro}
\proof
By the Cauchy-Riemann equations (\ref{CR}),
\begin{eqnarray}
\label{CRzw}
z_{\bar\zeta}=a \bar z_{\bar\zeta}, \qquad
w_{\bar\zeta}=b \bar z_{\bar\zeta}+c \bar w_{\bar\zeta}.
\end{eqnarray}
By the Chain Rule,
\begin{eqnarray*}
w_\zeta=w_z z_\zeta+w_{\bar z}\bar z_\zeta, \qquad
w_{\bar\zeta}=w_z z_{\bar\zeta}+w_{\bar z}\bar z_{\bar\zeta}.
\end{eqnarray*}
By eliminating
$z_{\bar\zeta}$,
$w_\zeta$, and
$w_{\bar\zeta}$, and canceling by $\bar z_{\bar\zeta}\ne0$
we obtain
\begin{eqnarray}
\label{wzwz}
w_{\bar z}-ac\bar w_z=-aw_z+c\bar w_{\bar z}+b.
\end{eqnarray}
Plugging (\ref{CR1}) in (\ref{wzwz}) yields
\begin{eqnarray}
\label{CR1relations}
a_1-ac\bar c_1=-a,\quad
b_1-ac\bar b_1=b,\quad
c_1-ac\bar a_1=c.
\end{eqnarray}
Solving (\ref{CR1relations}) yields the coefficients
of the desired equation (\ref{CR1}):
\begin{eqnarray}
\label{CR1coeff}
a_1=-a(1-|c|^2)/\Delta,\quad
b_1=(b+ac\bar b)/\Delta,\quad
c_1=c(1-|a|^2)/\Delta,\quad
\Delta=1-|ac|^2,
\end{eqnarray}
in which $|a_1|+|c_1|<1$.
Conversely, if $w(z)$ satisfies (\ref{CR1}),
then we can find the parameter $\zeta(z)$ by
solving the linear Beltrami equation
$\zeta_{\bar z}+a(z,w(z))\zeta_z=0$,
which is the first equation in (\ref{CRzw})
written for the inverse function $z\mapsto\zeta(z)$.

Finally, we observe that the equations (\ref{CR1relations})
have a unique solution in $(a,b,c)$ with $|a|<1, |c|<1$.
Indeed, put $q=ac$. Then (\ref{CR1relations}) imply
$(a_1-q\bar c_1)(c_1-q\bar a_1)=-q$. The latter
has only one solution $q$ with $|q|<1$ because $|a_1|+|c_1|<1$.
The equations (\ref{CR1relations}) now give the triple $(a,b,c)$
with $|a|<1, |c|<1$.
The other conclusions follow automatically.
Q.E.D.
\medskip

We need a maximum principle for solutions of (\ref{CR1}).
It is mentioned in \cite{Ve} without proof.
For completeness we include it here.
Consider a linear equation
\begin{eqnarray}
\label{CRq12}
w_{\bar z}=q_1 w_z + q_2 \bar w_{\bar z} + Q_1w+Q_2\bar w,
\end{eqnarray}
\begin{e-pro}
\label{Max-linear}
Let $q_1, q_2, Q_1, Q_2$ be bounded functions
in $\bar\D$ and let $|q_1|+|q_2|\le q_0<1$,
$|Q_1|+|Q_2|\le Q_0$, here $q_0$ and $Q_0$
are constants.
\begin{itemize}
\item[(i)]
There exists a constant $C>0$ depending only on $q_0$ and $Q_0$
so that for solutions of the equation (\ref{CRq12})
we have $\max_\D |w| \le C \max_{b\D} |w|$.
\item[(ii)]
If a solution to the equation (\ref{CRq12}) satisfies
$\Re w|_{b\D}=0$, then either $w$ does not vanish in $\bar\D$
or $w\equiv0$.
\end{itemize}
\end{e-pro}

We note that  $||Q_1||_{L^p}+||Q_2||_{L^p}\le Q_0$ for some
$p>2$ instead of $p=\infty$ would suffice, but we do not need
it here.
\proof
Let $w$ be a solution of (\ref{CRq12}).
Put
$q=q_1+q_2 \frac{\bar w_{\bar z}} {w_z}$ and
$Q=Q_1+Q_2 \frac{\bar w} {w}$.
Then $w$ also satisfies
\begin{eqnarray}
\label{CRq}
w_{\bar z}=qw_z + Qw.
\end{eqnarray}
Let $u$ and $v$ be solutions of the Beltrami
equations
\begin{eqnarray}
\label{Beltrami}
u_{\bar z}=qu_z,\qquad
v_{\bar z}=qv_z+Q
\end{eqnarray}
so that $u:\bar\D\to\bar\D$ is a homeomorphism.
Such solutions exist because $|q|\le q_0<1$.
Then one can see (\cite{Ve}, Theorem 3.31) that every solution
of (\ref{CRq}), whence (\ref{CRq12}), has a representation
\begin{eqnarray}
\label{CRq-rep}
w(z)=\phi(u(z))e^{v(z)},
\end{eqnarray}
where $\phi$ is holomorphic in $\D$.
Indeed, define $\phi$ by (\ref{CRq-rep}) and plug in (\ref{CRq}).
Then using (\ref{Beltrami}) we obtain
$(1-|q|^2)\bar u_{\bar z}e^v \phi_{\bar u}=0$.
Hence $\phi_{\bar u}=0$, and $\phi$ is holomorphic.

Although the equations (\ref{Beltrami}) depend
on a solution $w$, by Proposition \ref{Beltrami-Re0}(i)
there exists $v$ whose sup-norm
depends on $q_0$ and $Q_0$ only.
Hence, the assertion (i) follows by the usual maximum
principle applied to a holomorphic function $\phi$.

To prove (ii), we note that by Proposition \ref{Beltrami-Re0}(i),
the solution $v$ to the non-homogeneous equation
in (\ref{Beltrami}) can be
chosen with the condition $\Im v|_{b\D}=0$.
Then $\Re w|_{b\D}=0$ implies $\Re \phi|_{b\D}=0$.
Since $\phi$ is holomorphic, then $\phi\equiv ic$,
where $c\in\R$ is constant, and (ii) follows.
The proof is complete.
Q.E.D.

\begin{e-cor}
\label{Max-principle}
Let the coefficients $a_1$ and $c_1$ of (\ref{CR1})
satisfy $|a_1|+|c_1|\le a_0<1$ and
let $|(b_1)_w|+|(b_1)_{\bar w}|\le b_0$;
here $a_0$ and $b_0$ are constants.
Suppose $b_1(z,0)=0$.
Then there exists a constant $C>0$ depending
only on $a_0$ and $b_0$ so that for every solution
of (\ref{CR1}) we have
$
\max_\D |w| \le C \max_{b\D} |w|
$.
\end{e-cor}
\proof
Our hypotheses on $b_1$ imply that
$b_1(z,w)=Q_1(z,w)w+Q_2(z,w)\bar w$, where
$|Q_1|\le b_0$ and $|Q_2|\le b_0$.
Then the conclusion follows by Proposition \ref{Max-linear}(i).
Q.E.D.
\medskip

Finally we include a result on the regularity of the boundary
value problem for (\ref{CR1}).

\begin{e-pro}
\label{reg}
Let $w \in W^{1,p}(\D)$, $p > 2$, be a solution of (\ref{CR1})
satisfying the boundary condition $\Re w|_{b\D} = \phi$.
Suppose that $|a_1|+|c_1|< 1$,
$a_1,b_1,c_1 \in C^{k,\alpha}(\D \times \cc)$,
$0<\alpha<1$, and
$\phi \in C^{k+1,\alpha}(\D)$, $k \ge 0$.
Then $w \in C^{k+1,\alpha}(\D)$.
\end{e-pro}
\proof
We interpret the graphs of solutions of (\ref{CR1})
as $J$-complex discs and apply a reflection principle
from \cite{IvSu}. According to it, if $J$ is $C^{k,\alpha}$
($k\ge 0$, $0<\alpha<1$), then $J$-complex discs attached
to a totally real submanifold are $C^{k+1,\alpha}$-smooth.

Let $J$ be the almost complex structure corresponding
to the equation (\ref{CR1}) by Proposition \ref{CR1prop}.
Then the solution of (\ref{CR1}) with boundary
condition $\Re w|_{b\D} = \phi$ defines
a $J$-complex disc attached to the totally real
submanifold $\{ (z,w):|z|= 1, \Re w=\phi(z) \}$.
The conclusion now follows by the reflection
principle \cite{IvSu}.
Q.E.D.

\section{Proof of the main theorem}

Without loss of generality we can assume
$\Lambda=b\D\times b\D$.
Indeed, otherwise we can transform the curves
$\{z\}\times\gamma_z$ into the circles $\{z\}\times b\D$
by a smooth change of coordinates of the form
$z'=z$, $w'=w'(z,w)$ so that $w'(z,0)=0$ and
$w'=w$ for big $w$. By Proposition \ref{Change-abc},
the hypotheses of the theorem will be preserved by the change.

By Proposition \ref{CR1prop}, if a smooth $J$-complex disc
in $\bar\D\times\cc$
can be represented as a graph of a smooth function
$w:\bar\D\to\cc$, then $w$ satisfies the equation (\ref{CR1}).
The disc will be attached to $\Lambda=b\D\times b\D$
if in addition $|w(z)|=1$ for $z\in b\D$.

Without loss of generality we can assume that
the coefficients $a$, $b$, and $c$ have compact support
in $\bar\D\times\cc$. Indeed, by the maximum principle
(Corollary \ref{Max-principle}),
the boundary condition
$|w|=1$ implies $|w|<C$ for all $z\in\D$.
Here $C>0$ depends only on the coefficients
of (\ref{CR1}). Then we can assume $a=b=c=0$,
hence $a_1=b_1=c_1=0$ for $|w|>C$.

Since we are looking for a non-vanishing solution $w$,
we now make a change of variable $w=e^u$.
Then the equation (\ref{CR1}) with boundary
condition $|w|=1$ transforms into
\begin{eqnarray}
\label{CR2}
u_{\bar z}=a_2 u_z + c_2 \bar u_{\bar z} + b_2
\end{eqnarray}
with boundary condition
\begin{eqnarray}
\label{boundary}
\Re u(z)=0, z\in b\D.
\end{eqnarray}
Here
$$
a_2(z,u)=a_1(z,e^u),\quad
b_2(z,u)=e^{-u}b_1(z,e^u),\quad
c_2(z,u)=e^{\bar u-u}c_1(z,e^u).
$$
The coefficients of (\ref{CR2}) still
have uniformly bounded derivatives in $u$ and satisfy the
ellipticity condition $|a_2|+|c_2|<a_0<1$.
Under these hypotheses by \cite{Mo} (pp. 335--351;
see also \cite{SuTu1}, Proposition 4.2)
for every every $z_0\in b\D$ and $\tau\in\R$,
the boundary value problem (\ref{CR2}-\ref{boundary})
has a unique solution with $u(z_0)=i\tau$.

Returning to the original variable $w$, we obtain
a one-parameter family of discs attached to $\Lambda$.
They are parametrized by $\tau\in\R/2\pi\Z$ for fixed $z_0$,
say $z_0=1$.
Their boundaries are disjoint and cover all of $\Lambda$.
Since $w=e^u\ne 0$, then they do not meet $\bar\D\times\{0\}$.
By the continuous dependence statement of \cite{Mo},
the family is continuous. For the same reason, they
depend continuously on $A$, which gives the continuous
dependence conclusion in (ii).

The smoothness of the above one-parameter family follows by
the implicit function theorem. Fix $k \geq 1$, $0<\alpha<1$.
By Proposition \ref{reg}, every solution $u\in C^{k+1,\alpha}(\D)$.
Denote by $L: C^{k,\alpha}(\D) \to C^{k-1,\alpha}(\D)$,
$L: \dot u \mapsto  L(\dot u)$, the linearization at $u$
of the operator defined by (\ref{CR2}).
Then the operator $L$ itself can be written in the same form
with coefficients of class $C^{k-1,\alpha}(\D)$.
By Proposition \ref{reg} the linear map
${\cal L}:C^{k,\alpha}(\D)\to C^{k-1,\alpha}(\D)
\times C^{k,\alpha}(b\D)$ defined by
${\cal L}(\dot u) = (L(\dot u), \Re \dot u \vert_{b\D})$
is surjective. As it was shown above, the solution of the boundary
value problem (\ref{CR2}-\ref{boundary}) is uniquely
determined by the condition $u(z_0) = i\tau$ for
every $z_0\in b\D$ and $\tau\in\R$. Hence the kernel
of ${\cal L}$ is one-dimensional. By the implicit function
theorem there exists a $C^{k,\alpha}$-smooth one-parameter
family of solutions of (\ref{CR2}-\ref{boundary}).
By uniqueness we obtain the smooth dependence
conclusion in (ii).

The fact that the discs are disjoint
follows by the positivity of intersections of
$J$-complex curves. Indeed, we can include our
boundary value problem $|w|=1$ for the equation
(\ref{CR1}) in a one-parameter family with
boundary condition $|w|=r$.
Let $f$ be a $J$-complex disc constructed above,
attached to $\Lambda$.
Then we can construct a continuous family
of $J$-complex discs $f_r$ with boundary condition $|w|=r$,
so that $f_0$ is the disc $\bar\D\times\{0\}$
and $f_1=f$. Let $g$ be another disc constructed above,
attached to $\Lambda$. Since the boundaries
of $g$ and $f_r$ do not intersect, then the
intersection index of $f_r$ and $g$ is independent
of $r$. Since $g$ does not intersect $f_0$,
then it does not intersect $f_1=f$ either.

We now prove that the surface $\Gamma$ swept out by the
discs is smooth. Let $z\mapsto u(z,\tau)$ be the solution
of (\ref{CR2}) with conditions $\Re u|_{b\D}=0$ and
$u(1,\tau)=i\tau$.
It suffices to show that the map
$(z,\tau)\mapsto (z,u(z,\tau))$ is an immersion, which
in turn reduces to $\partial u/\partial\tau\ne0$.
Put $v=\partial u/\partial\tau$ and differentiate
(\ref{CR2}) with respect to $\tau$. Then we obtain
a linear equation
\begin{eqnarray}
\label{CR3}
v_{\bar z}=a_2 v_z + c_2 \bar v_{\bar z}
+Q_1v+Q_2\bar v
\end{eqnarray}
with boundary condition $\Re v|_{b\D}=0$.
Here
$Q_1=(a_2)_u u_z + (c_2)_u \bar u_{\bar z} + (b_2)_u$ and
$Q_2=(a_2)_{\bar u} u_z + (c_2)_{\bar u} \bar u_{\bar z}
+(b_2)_{\bar u}$.
Since $v(1)=i\ne0$,
then by Proposition \ref{Max-linear}(ii),
the solution $v$ does not vanish in $\bar\D$ as desired.

Finally, the part (iii) follows by Corollary \ref{Max-principle},
the maximum principle for (\ref{CR1}). The proof is complete.
Q.E.D.

\end{document}